\documentclass{ifacconf}

\usepackage{graphicx}      
\usepackage{natbib}        

\usepackage{amssymb,amsmath,bm}
\usepackage[dvipsnames]{xcolor}

\usepackage{graphicx,float}
\newcommand{\beq}{\begin{equation}}
\newcommand{\eeq}{\end{equation}}

\usepackage{algorithm}
\usepackage{algpseudocode}
\newcommand{\lu}{\rule{0.15 cm}{.5pt}} 
\DeclareMathAlphabet{\mathcalo}{OMS}{cmsy}{m}{n} 
\begin{document}
\begin{frontmatter}
This work has been submitted to IFAC for possible publication
\title{Sampling-Based Risk-Aware Path Planning Around Dynamic Engagement Zones} 

\thanks[footnoteinfo]{This work was supported by the AFOSR Summer Faculty Fellowship Program and AFOSR 24RQCOR002. \\ Distribution Statement A.  Approved for public release: distribution is unlimited. AFRL-2024-0898 , Cleared 20 FEB 2024.}

\author[First]{Artur Wolek} 
\author[Second]{Isaac E. Weintraub} 
\author[Second]{Alexander Von Moll} 
\author[Second]{David Casbeer} 
\author[Third]{Satyanarayana G. Manyam}

\address[First]{Department of Mechanical Engineering and Engineering Science, \\ University of North Carolina at Charlotte, Charlotte, NC, 28223 USA (e-mail: {awolek@charlotte.edu}).}
\address[Second]{Air Force Research Laboratory, Wright-Patterson AFB, OH 45433, USA  (e-mail: {\{isaac.weintraub.1, alexander.von\_moll, david.casbeer\}@us.af.mil} )}
\address[Third]{DCS Corp., 4027 Col Glenn Hwy, Dayton, OH, 45431, USA (e-mail: {msngupta@gmail.com})}

\begin{abstract}
Existing methods for avoiding dynamic engagement zones (EZs) and minimizing risk leverage the calculus of variations to obtain optimal paths. While such methods are deterministic, they scale poorly as the number of engagement zones increases. Furthermore, optimal-control based strategies are sensitive to initial guesses and often converge to local, rather than global, minima. This paper presents a novel sampling-based approach to obtain a feasible flight plan for a Dubins vehicle to reach a desired location in a bounded operating region in the presence of a large number of engagement zones. The dynamic EZs are coupled to the vehicle dynamics through its heading angle. Thus, the dynamic two-dimensional obstacles in the $(x,y)$ plane can be transformed into three-dimensional static obstacles in a lifted $(x,y,\psi)$ space. This insight is leveraged in the formulation of a Rapidly-exploring Random Tree (RRT$^*$) algorithm. The  algorithm is evaluated with a Monte Carlo experiment that randomizes EZ locations  to characterize the success rate and average path length as a function of the number of EZs and as the computation time made available to the planner is increased.
\end{abstract}
\begin{keyword}
risk-aware path planning, dynamic obstacles, rapidly-exploring random trees
\end{keyword}

\end{frontmatter}

\section{Introduction}
Air and missile threats in contested airspace represent a challenge for a number of Air Force missions \citep{scott2017countering, odonohu2019joint}. A common objective for navigating in threat-laden environments is to reach some desired location in minimum time while minimizing risk of capture. This work proposes an approach for finding minimum time paths for a turn-constrained vehicle through a field of dynamically changing engagement zones (EZs).

Guided threats reach their intended targets in a dynamic fashion; a number of guidance laws exist and can be found in the literature, for example, see \citep{zarchan2012tactical}. Many of these systems leverage the line-of-sight angle to attain capture \citep{zarchan1999ballistic}. The reachable set of these systems, defining the EZ, may be obtained by propagating their equations of motion. In this work, a cardioid is leveraged as a surrogate EZ since its shape is representative  \citep{herrmann1996air, yan2020reachability}. Tuning the shape to match an existing threat model is out of the scope of this work; but, doing so would inform mission planning around threat regions.

In this paper, engagement zones are dynamic in nature and represent the 
the region of space whereby if the aircraft were to keep course it would be in threat of capture by a designated pursuer. As the aircraft navigates to its intended terminal location, the EZs change dynamically as a function of the aircraft's relative pose. A discussion of various EZ models when the aircraft exhibits simple motion is found in \citep{von2023basic}. 

Prior work focusing on penetrating EZs has considered minimizing detection from radar sites using the A* algorithm \citep{zhang2020novel, zhang2023efficient} and in 3-D using particle swarm optimization, grey wolf optimization, and genetic algorithms \citep{xu2020optimized}. Other approaches that account for the dynamic nature of a single EZ \citep{weintraub2022optimal} and two EZs \citep{dillon2023optimal} have leveraged the calculus of variations and direct optimal control. These approaches do not scale well as the number of EZs increases, suffer from sensitivity due to initial guesses, and typically resolve to a locally optimal solution.

\begin{figure*}
\centering
\includegraphics[width=0.8\textwidth]{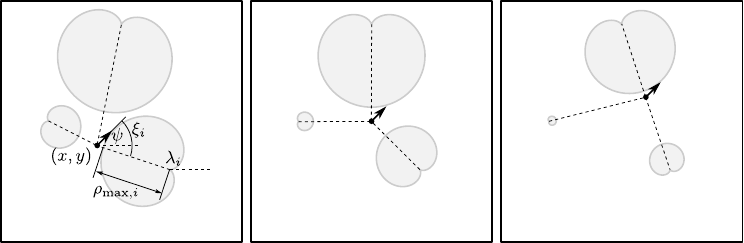}
\caption{Example scenario with three dynamic EZs and an aircraft moving along a straight path. Each EZ is drawn by sweeping out $\theta_i \in [0, 2\pi)$ in \eqref{eq:rho_Rmin_zero}. The EZ geometry depends on the relative $\lambda_i$ and $\xi_i$ that vary with aircraft motion. }
\label{fig:setup}
\end{figure*}

The Rapidly-Exploring Random Tree (RRT) motion planning algorithm  was first introduced by \cite{lavalle2001randomized}. 
The main idea of this algorithm is to incrementally build a tree rooted at the start node. Samples are randomly generated in the free space and connected to the tree whenever possible. The tree is grown until a feasible sample is generated in the goal region such that it can be connected to the tree. \cite{karaman2011sampling} later presented RRT$^*$ and proved that it finds the optimal solution asymptotically.  Another variant, informed RRT$^*$ \citep{gammell2014informed}, generates random samples from a subset of the free-space to improve efficiency while also preserving asymptotic optimality. The FMT$^*$ (Fast Marching Tree) algorithm in \citep{janson2015fast} combines the ideas of RRT$^*$ and probabilistic road-maps. The algorithm converges faster by reducing the number of collision checks. The algorithm presented in \citep{gammell2020batch}, BIT$^*$ combines the benefits of graph search and sampling-based methods. The algorithms above consider static obstacles, and motion planning problems with dynamic obstacles are considered by \cite{otte2016rrtx}. They present a single query algorithm (no precomputing needed), RRT$^X$, which is asymptotically optimal. In a related work, path planning around static obstacles using a Dubins vehicle was presented in \citep{boardman2018improved}. The work presented in this paper considers dynamic EZs by lifting the space into an additional dimension: aircraft heading.

The contributions of this paper are: (1) identifying a lifted space for which a class of dynamic engagement zones (EZs) can be viewed as static obstacles, and (2) a sampling-based motion planning algorithm to plan a risk-aware path to avoid a large number of EZs. The performance of the algorithm is characterized through a Monte Carlo study that plans trajectories over randomized EZ locations.

The remainder of the paper is organized as follows. Section~\ref{sec:prob_formulation} formulates the risk-aware path planning problem.  Section~\ref{sec:rrt} describes the proposed RRT$^*$ sampling-based approach.  Section~\ref{sec:monte_carlo} presents Monte Carlo simulation results. Section~\ref{sec:conclusion} concludes the paper.

\section{Problem Formulation}
\label{sec:prob_formulation}
\subsection{Vehicle and Engagement Zone Models}
\label{sec:dynEZ}
Consider a constant-altitude fixed-wing aircraft modeled according to \citep{dubins1957curves}:
\begin{align}
\dot x(t) &= v \cos \psi(t) \;,\label{eq:xdot1} \\
\dot y(t) &= v \sin \psi(t) \;, \\
\dot \psi(t) &= u(t) \label{eq:thdot} \;, 
\end{align}
where $(x,y)$ is the planar position, $\psi$ is the heading, $v$ is a constant speed and $u$ is a bounded turn-rate, $|u| \leq u_{\rm max}$. 
The aircraft operates in an operating region $\mathcal{D} = [0,1]^2$ in the presence of $N$ dynamic engagement zones (EZs) located at planar positions $(x_i, y_i) \in \mathcal{D} $ for $i=1,\ldots,N$. 
This paper adopts the geometry of a dynamic EZ  proposed in \citep{weintraub2022optimal} with the outer boundary of the EZ (i.e., distance from center point $(x_i,y_i)$) described by a cardioid. In polar coordinates, the radius of the cardioid is given by
\begin{align} 
\rho (\theta_i, \lambda_i, \xi_i) &= \left( 
\left(
\frac{\cos \xi_i + 1}{2}
\right)
(R_{\rm max} - R_{\rm min}) + R_{\rm min}
\right) \notag \\
& \qquad 
\frac{1}{2}
\left(
1 + \sin
\left(
\frac{\pi}{2} - \lambda_i + \theta
\right)
\right) \;,
\label{eq:rho_th_lam_xi}
\end{align}
where $\theta \in [0, 2\pi)$ is the polar angle, 
\begin{align}
\xi_i (x,y,\psi) = \psi - \lambda_i  (x,y)- \pi\;,
\label{eq: relative bearing}
\end{align}
 is the relative bearing angle, $R_{\rm max}$ is the maximum radius of the EZ, $R_{\rm min}$ is the minimum radius of the EZ, and 
\begin{equation}
    \lambda_i(x,y) = {\rm atan}((y-y_i)/(x-x_i))\;,
    \label{eq:lambda}
\end{equation}
is the line-of-sight angle from the center of tbe $i^{\rm th}$ EZ to the aircraft.  Following  \cite{weintraub2022optimal}, the minimum range is selected as $R_{\rm min} =0$ so that \eqref{eq:rho_th_lam_xi} reduces to:
\begin{align} 
\rho (\theta_i, \lambda_i, \xi_i) = 
\tfrac{R_{\rm max}}{4}
\left( 
\cos \xi_i \hspace{-0.1em} + \hspace{-0.1em} 1
\right)
\left(
1 \hspace{-0.1em} + \hspace{-0.1em} \sin
\left(
\tfrac{\pi}{2} - \lambda_i + \theta
\right)
\right).
\label{eq:rho_Rmin_zero}
\end{align}
The only point along the curve \eqref{eq:rho_th_lam_xi} that is potentially intersecting the aircraft is the point where the polar angle is equal to the line-of-sight (i.e., $\theta_i(t) = \lambda_i(t)$). With these assumptions \eqref{eq:rho_Rmin_zero} becomes \citep{weintraub2022optimal}:
\begin{align} 
\rho_{{\rm max},i}(\xi_i) 
&= \tfrac{R_{\rm max}}{2}( \cos \xi_i + 1) \;,
\label{eq:rho_max_xi}
\end{align}
which describes the worst-case (maximum range) outer boundary of a cardioid shape as a function of relative bearing (see Fig.~\ref{fig:setup}). The distance of the aircraft from the center of the $i^{\rm th}$ EZ is given by $d_i(x,y) = \sqrt{(x-x_i)^2 + (y-y_i)^2}$.
\begin{figure*}[htb]
\centering
\includegraphics[width=0.9\textwidth]{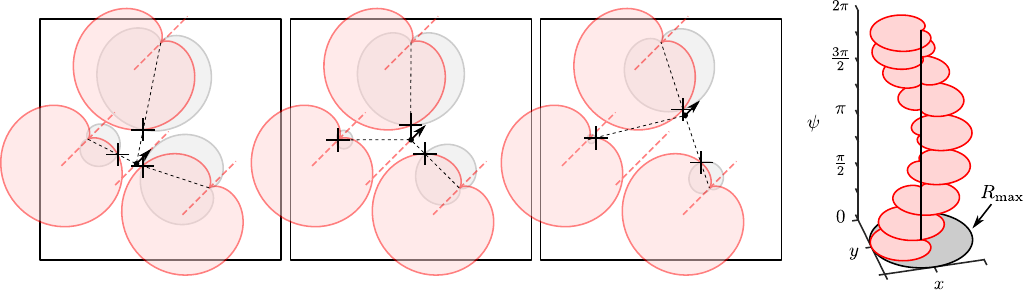}
\caption{Left three panels: Example scenario (identical to Fig.~\ref{fig:setup}) with three EZ regions and aircraft  moving along a straight line path. The red colored EZ obstacles are drawn by varying $\lambda_i \in [0, 2\pi)$ in \eqref{eq:rho_lam_psi} for a fixed heading angle $\psi = \pi/4$. The ``+" markers indicate locations at which both the EZ obstacle and dynamic EZ intersect. These points correspond to locations that are a distance $\rho_{{\rm max},i}$ from the aircraft. Right panel: Slices of the EZ obstacle at different heading planes in the lifted space.}
\label{fig:setup_obs}
\end{figure*}
\subsection{Problem Statement}
The problem investigated in this paper is to plan a minimum-time path for an aircraft to join an initial and terminal state while remaining in the operating region and avoiding the dynamic EZs. That is, the optimal control problem is to
\begin{equation}
    ~{\rm minimize}\qquad J({u}(\cdot)) = \int_0^T 1 \; \mathrm{d}t = T
\end{equation}
where $J$ is the cost functional, $T$ is the final time, and $u(\cdot)$ is a piecewise-continuous turn-rate control, subject to the kinematic constraints \eqref{eq:xdot1}--\eqref{eq:thdot}, the boundary conditions
\begin{equation}
\begin{aligned}
x(0) &= x_0   \\
y(0) &= y_0\\
\psi (0) &= \psi_0 \\
\end{aligned}
\qquad  \qquad \qquad 
\begin{aligned}
x(T) &= x_T   \\
y(T) &= y_T\\
\psi (T) &= \psi_T \\
\end{aligned} \;,
\end{equation}
and the state constraints
\begin{align}
( x(t), y(t) )^\top &\in \mathcal{D} \label{eq:constraint0}\\
g_i(x,y,\psi) := \rho_{\rm max}(\xi(x,y,\psi)) - d_i(x,y) &\leq 0 \label{eq:constraint1}
\end{align}
for all $i = 1, \ldots, N$, 
 and the control constraint $|u(t)| \leq u_{\rm max}$ for all time $t\in[0,T]$.

\section{Risk-Aware RRT$^*$  Path Planning Algorithm}
\label{sec:rrt}
\subsection{Engagement Zone Obstacles}
At first glance, the constraint \eqref{eq:constraint1} appears to be a path constraint that depends non-linearly on the state variables. In the $(x,y)$ plane the EZ region scales and rotates as the aircraft changes pose (see Fig.~\ref{fig:setup}). The typical way to handle path constraints in an optimal control problem is using the method of Lagrange multipliers.  However, the constraint $g(x,y,\psi) = 0$ can alternatively be viewed as a set of  \emph{static} obstacles (as opposed to the dynamic EZ described in Section \ref{sec:dynEZ}) in the three-dimensional ``lifted" $(x,y,\psi)$ space. (Note, this space is  the  configuration space of the Dubins vehicle, $\mathcal{Q} = \mathcal{D} \times [0, 2\pi)$.) Related methods have been used to compute  configuration space obstacles that consider the shape and size of a mobile robot in the presence of static obstacles using Minkowski sums  \citep{lozano1983spatial}. Our approach is similar, but the mobile robot is a point and obstacle geometries are heading-dependent. The configuration-space perspective with static obstacles motivates the use of path planning methods that use collision checking routines against geometric obstacles. 

To construct the EZ obstacle, substitute \eqref{eq: relative bearing} into \eqref{eq:rho_max_xi} and treat the heading as a constant parameter: 
\begin{align}
    \rho_{{\rm max},i}(\lambda_i; \psi) &= \tfrac{R_{\rm max}}{2} ({1-\cos(\psi - \lambda_i) }) \;.
\label{eq:rho_lam_psi}
\end{align} 
The cardiod generated by sweeping out $\lambda_i \in [0, 2\pi)$ in \eqref{eq:rho_lam_psi}, for a fixed value of $\psi$, is different from the cardiod generated by sweeping $\psi \in [0, 2\pi)$ in \eqref{eq:rho_lam_psi}, for a fixed value of $\lambda_i$. This new cardiod is a cross-section of the \emph{EZ obstacle} on a particular heading plane $\psi$. The three-dimensional EZ obstacle in the configuration space $\mathcal{Q}$ is visualized in the right panel of Fig.~\ref{fig:setup_obs} and is defined as the union of such  cardiods \eqref{eq:rho_lam_psi} across all heading planes $\psi \in [0, 2\pi)$. The set of all EZ obstacles is:
\begin{align}
\mathcal{O} &= \{ (x,y,\psi) \in \mathcal{Q}~|~d_i(x,y) \leq  \rho_{{\rm max},i}(\lambda_i(x,y), \psi)   \},
\label{eq:ez_obstacles}
\end{align}
for $i = 1, \ldots, N$.
On each constant heading plane the cardiod has the same size and shape and is symmetric about the angle $\psi$ drawn on that plane; however, as the heading changes the EZ rotates around the vertical axis passing through $(x_i, y_i)$. If the EZ obstacle is projected onto the $(x,y)$ plane it gives a circle of radius $R_{\rm max}$ (i.e., a conservative bound on the geometry of the EZ). 

\subsection{Dubins Motion in the Configuration Space} 
\begin{figure}[htb]
\centering
    \includegraphics[width = 0.47\textwidth]{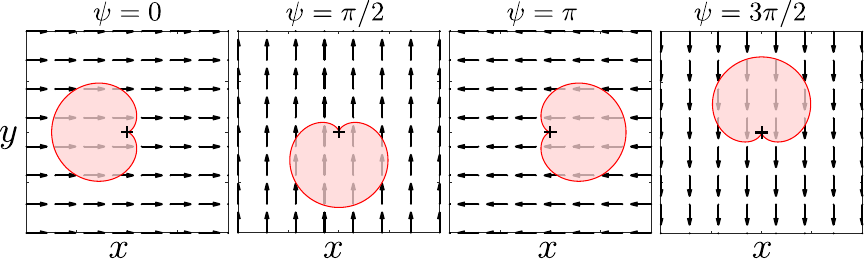}
    \caption{The engagement zone rotates on each heading plane (indicated by the vector field of allowable motion).}
\label{fig:four_headings}
\end{figure}

In the configuration space, maximum turn-rate arcs for the Dubins model map to helices that have a radius of $R = v/u_{\rm max}$. During one full revolution of the helix, the vehicle travels from the 0 heading plane (bottom) to the $2\pi$ heading plane (top). In other words, the helices have a pitch of $2\pi$. All straight segments are along constant heading planes and, on each slice of fixed heading $\psi$, the Dubins vehicle can only move on straight lines in the direction $\psi$, as illustrated by the vector fields in Fig.~\ref{fig:four_headings}. Dubins paths, therefore, appear as helix-straight-helix segments or helix-helix-helix segments. The transition between segments are piece-wise continuous with corners in-between segments when viewed in the configuration space. The corners appear since the curvature is discontinuous at the junctions of turns in a Dubins path. 

The configuration space offers a more intuitive (or ``lifted") view of the problem. For example, simply plotting the obstacles in the 3D configuration space can reveal if a feasible path does not exist. That is, if EZ obstacles intersect to form a gapless boundary around the vehicle then there is no feasible solution. However, if the EZ obstacles do contain a gap we cannot conclude that a feasible Dubins path exists. The Dubins vehicle has turning dynamics that may limit it from maneuvering into a position to pass through this gap. 

\subsection{Risk-Aware RRT$^*$ Path Planning Algorithm}
The  Rapidly-exploring Random Tree (RRT$^*$) algorithm is a sampling-based algorithm that is well-suited towards planning in obstacle-rich environments.  
Below, we briefly summarize the RRT$^*$ algorithm (following \cite{karaman2011sampling}) using the pseudocode in Alg.~\ref{alg:RRT_pseudo} and highlight how it relates to our problem. In practice, we adopt the RRT$^*$ implementation in C++ using the Open Motion Planning Library (OMPL) \citep{sucan2012open} which may differ slightly from  Alg.~\ref{alg:RRT_pseudo}.

RRT$^*$ constructs a tree $G = (V,E)$  with nodes $V = \{ {\bm q}_0, {\bm q}_1, \ldots, {\bm q}_{|V|}\}$ that represent aircraft states (i.e., ${\bm q}_i \in \mathcal{Q}$) and with directed edges  $E \subseteq V^2$  (denoted by pairs of states) that  represent feasible paths. For example, the edge $({\bm q}_i,{\bm q}_j) \in E$ denotes a feasible path from ${\bm q}_i$ to ${\bm q}_j$ that satisfies the dynamics \eqref{eq:xdot1}--\eqref{eq:thdot} and constraints \eqref{eq:constraint0}--\eqref{eq:constraint1}. The cost of an edge is the associated travel time.

The tree is initialized with a single root node equal to the initial state ${\bm q}_0 = (x_0,y_0,\psi_0)^\top$ (Alg.~\ref{alg:RRT_pseudo}, line 2) and, while the elapsed computation time $\tau$ is less than a user-defined limit $\tau_{\rm max}$ (line 3), the algorithm incrementally adds nodes and edges to the tree to reach the goal ${\bm q}_{\rm goal} = (x_T, y_T, \psi_T)^\top$. Even after the goal is reached, RRT$^*$ continues to evaluate alternative, possibly lower-cost,  solutions. First, the SAMPLE\lu FREE function (line 4) randomly selects a state ${\bm q}_{\rm rand}$ in the free space (i.e., that is in the operating region and outside the EZ obstacles \eqref{eq:ez_obstacles}). Next, the state ${\bm q}_{\rm nearest}$  in the set $V$ that is nearest to  ${\bm q}_{\rm rand}$ (according to Euclidean distance) is identified (line 5). Then, a Dubins path is planned from ${\bm q}_{\rm near}$ to ${\bm q}_{\rm rand}$ using the STEER subroutine. This subroutine returns a point ${\bm q}_{\rm new}$ that is partway along the planned Dubins path (e.g., after a fixed arclength stepsize). The EZ engagement is handled by the OBSTACLE\lu FREE subroutine (line 7) that checks if the segment of the path from ${\bm q}_{\rm nearest}$ to ${\bm q}_{\rm new}$ is in collision on the appropriate heading plane and that the boundary of the operation region is not violated.  
NEAR (line 8) finds a set of nodes $Q_{\rm near}$ in the graph that are in a neighborhood around ${\bm q}_{\rm new}$ according to a criteria such as relative Euclidean distance. The OPT\lu PARENT function then determines the parent node ${\bm q}_{\rm min}$ in the set $Q_{\rm near} \cup {\bm q}_{\rm nearest}$ that gives the lowest cost total path (from the starting location) to ${\bm q}_{\rm new}$ (line 9). The new node and new edge are added to the graph (lines 10, 11). Lastly, the REWIRE step (line 12) searches if any node ${\bm q} \in Q_{\rm near}$ can be reached with lower cost by utilizing the new node and edge added --- if so, the tree is re-wired by adding and removing edges appropriately. 

\begin{algorithm}[h!]
    \caption{Risk-Aware RRT$^*$ Planner}
   \label{alg:RRT_pseudo}
    \begin{algorithmic}[1] 
        \Procedure{Risk-Aware RRT$^*$}{${\bm q}_0$, ${\bm q}_{\rm goal}$, $\tau_{\rm max}$}  
            \State $G$.init(${\bm q}_0$) 
            \While{$\tau$ $\leq \tau_{\rm max}$}
                \State ${\bm q}_{\rm rand} \gets$ SAMPLE\lu FREE$( \mathcal{O} )$
                \State ${\bm q}_{\rm nearest} \gets$ NEAREST\lu CONFIG$(G, {\bm q}_{\rm rand})$ 
                \State ${\bm q}_{\rm new} \gets$ STEER$({\bm q}_{\rm nearest}, {\bm q}_{\rm rand})$  
        	       \If { OBSTACLE\lu FREE(${\bm q}_{\rm nearest},{\bm q}_{\rm new}$, $\mathcal{O}$) }  
        	         \State $Q_{\rm near} \gets$ NEAR$(G, {\bm q}_{\rm new})$
                  \State ${\bm q}_{\rm min} \gets$ OPT\lu PARENT$(G, Q_{\rm near}, {\bm q}_{\rm new})$
                  \State $G$.add\lu node(${\bm q}_{\rm new}$)  
        	         \State $G$.add\lu edge(${\bm q}_{\rm min}, {\bm q}_{\rm new}$) 
                  \State $G \gets$ REWIRE$(G, Q_{\rm near}, {\bm q}_{\rm rand})$ 
        	       \EndIf
                
            \EndWhile\label{euclidendwhile}
            \State \textbf{return}  
        \EndProcedure
    \end{algorithmic}
\end{algorithm}

\subsection{Illustrative Example}
\label{sec:example}

This section presents an example of the proposed approach. The vehicle initial state is ${\bm q}_0 =  (0,0,0)^\top$ and the desired terminal state at the free-final time, $T$, is ${\bm q}_{\rm goal} = (1,1,0)^\top$. The length unit (LU) and time unit (TU) are normalized such that the vehicle's speed is $v = 1$ LU/TU. The turn radius of the Dubins vehicle is $R=0.1$ LU. A total of $N=16$ EZ regions were randomly generated in the operating domain. The proposed method does not require the EZ ranges to all be homogeneous, but for simplicity we assume that the maximum EZ range is $R_{\rm max} = 0.15$ for all EZs. (Also, while not presented here, the proposed methodology can be extended to other EZ geometries that exhibit aircraft heading-dependence, such as those outlined in \citep{von2023basic}.)

\begin{figure}[hhtb]
\centering
\includegraphics[width = 0.48\textwidth]{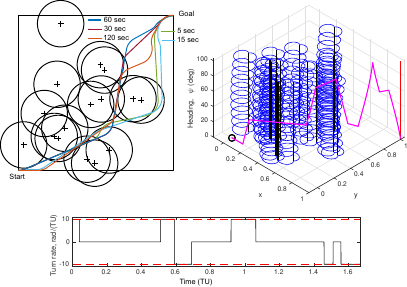}
    \caption{Left: Example of RRT solution with superimposed circles of $R_{\rm max}$ at WEZ locations. Right: The same trajectory in the configuration space $(x,y,\psi)$ with static EZ obstacles. Bottom: Corresponding turn-rate.}
\label{fig:example}
\end{figure}

The geometry of the randomly placed EZs is shown in the top left panel of Fig.~\ref{fig:example}. This panel draws a disc of radius $R_{\rm max}$ around each EZ point. This disc corresponds to a conservative bound on the dynamic EZ obstacle.  If one were to use these conservative obstacles, it would appear there is no solution to pass through the EZs. However, the solutions computed by the proposed method (shown for a range of computation budgets $\tau_{\rm solve}$) pass through the $R_{\rm max}$ circles. As the available computation time increases, the solution improves and changes topology based on which side (left or right) it passes  the right-most EZs. 

The top right panel of Fig.~\ref{fig:example} illustrates the 120 second solution in the lifted space with the static EZ obstacles. The trajectory consists of helices and straight line segments and the turns appear as sharp corners  due to the discontinuous turn-rate of Dubins curves. (However, due to the choice of relatively small turn radius the helical nature of the turning segments is not pronounced in Fig.~\ref{fig:example}.) In this paper we focused on the problem of a specific terminal heading. A more general problem could consider reaching the terminal point with a free heading angle (represented by the vertical red line in the top right panel of  Fig.~\ref{fig:example}). 

\begin{figure}[htb]
\centering
 \includegraphics[width = 0.40\textwidth]{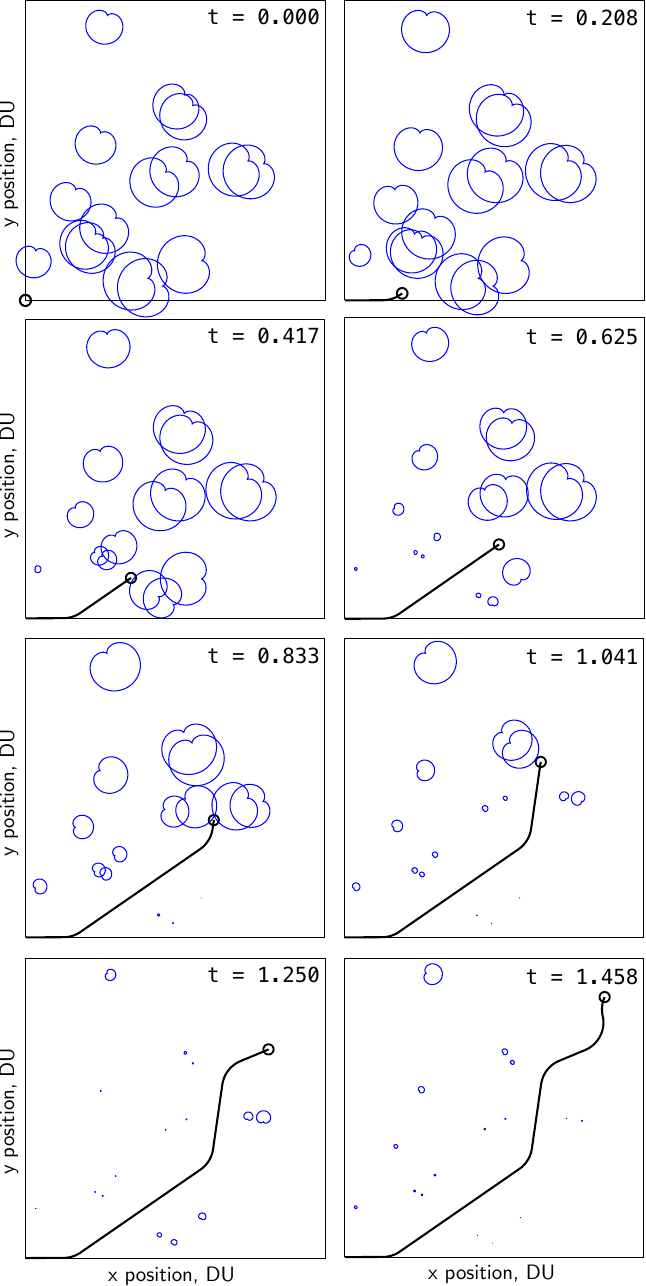}
    \caption{Eight snapshots in time of the aircraft traversing the path  in Fig.~\ref{fig:example} and the corresponding dynamic EZs. The aircraft position  at each timestamp shown in the upper right corner is denoted by a circular marker. }
\label{fig:example_snapshots}
\end{figure}

The bottom panel of Fig.~\ref{fig:example} provides the turn-rate control corresponding to the 120 second solution. The vehicle makes six turns along its path. Lastly, Fig.~\ref{fig:example_snapshots} illustrates several snapshots in time of the 120 second solution. This is the perspective taken in prior work \citep{weintraub2022optimal, dillon2023optimal} where the EZs are dynamic obstacles in the plane that scale and rotate.

\section{Performance Analysis}
\label{sec:monte_carlo}
To analyze the performance of the proposed approach, a Monte Carlo experiment was conducted that varied the number of EZs and available solve time for the RRT$^*$ algorithm. For each number of EZs, $N \in \{ 4, 8, 20, 24 \}$, a total of $M = 500$ different scenarios were generated. Half the EZs in each scenario were placed uniformly at random in the operating region and the other half were place according to a 2D Gaussian probability density function with a standard deviation of 1/5. (This mixed randomization strategy was adopted to make a more challenging set of scenarios by concentrating more EZs at the center of the operating region.) However, some random EZ placements may lead to ill-posed problems. For example, if the EZs are configured in such a way that they block a passage to the terminal state or if the vehicle is initialized inside an EZ. 
To ensure that scenarios generated were feasible,  only the EZ  configurations that were successfully solved with the solver using a $\tau_{\rm solve} =$ 160 second computation budget are used. The vehicle speed, EZ parameter $R_{\rm max}$, and the initial/terminal configurations were the same as those described in the illustrative example of Sec.~\ref{sec:example}. That is, the vehicle travels from the bottom left corner to the top right corner of the operating domain. The solver was run on each trial with a series of computation budgets $\tau_{\rm solve} \in \{ 5, 10, 20, 40, 60, 80, 100, 120, 140, 160\}$ seconds. 

After running the algorithm with each $(N, \tau_{\rm solve})$ setting, the total number of successfully solved problems (out of $M$) was recorded along with the reported cost (i.e., total transit time or, equivalently, path length). The solutions were verified by ensuring they satisfy the dynamics and do not intersect any obstacles or the boundary of the domain. 

\begin{figure}[htb]
\centering
\includegraphics[width=0.35\textwidth]{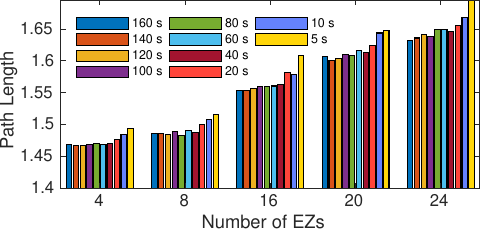}
\includegraphics[width=0.35\textwidth]{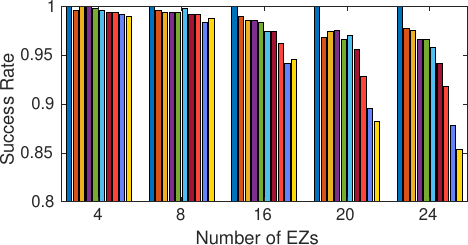}
\caption{Monte Carlo results path length (top) and success rate (bottom) as a function of maximum allowable RRT solve time $\tau_{\rm solve}$ with increasing number of EZs. }
\label{fig:monte_carlo}
\end{figure}

The top panel of Fig.~\ref{fig:monte_carlo} shows the mean path length increasing with the number of EZs. As expected, providing a higher computation time budgets allows leads to reduced costs. Further, as the number of EZs increases more maneuvering is required and the overall average path length increases. 
The bottom panel of Fig.~\ref{fig:monte_carlo} shows the success rate of the algorithm with increasing number of EZs. The success rate for the solver with $\tau_{\rm solve} = 160$~s is 100\% since a criteria for inclusion was that $\tau_{\rm solve} = 160$~s was feasible. A more informative result is given for $\tau_{\rm solve} = 120$~s --- the success rate varied from 96.8 \% to 99.6 \%. On the other extreme, the success rate for the solver with $\tau_{\rm solve} = 5$~s ranges from 85.4\% to 99.0\%. 

The results suggest that if computation time is critical (e.g., solutions are required in 20 seconds or less), then the RRT$^*$ algorithm proposed is only suitable for a small number of EZs (e.g., up to eight). For applications that can tolerate larger planning times the RRT$^*$ may be viable for a larger number of EZs. However, the RRT$^*$ only provides a  probabilistic completeness guarantee and thus it is possible that it does not return a solution, even when one exists, regardless of computation time budget. 

\section{Conclusion}
\label{sec:conclusion}
This paper investigated the problem of planning a risk-aware path through a large number of dynamic engagement zones with a Dubins vehicle model. It was recognized that the path constraints imposed by the EZs can be viewed as static obstacles in the configuration space, and that existing sampling-based motion planners with collision checking can be applied as an alternative to optimal-control-based solutions previously investigated. A  RRT$^*$ algorithm was implemented using the OMPL software package. The algorithm quickly finds a feasible solution and then continuously refines it for a user-specified time. The method was evaluated through Monte Carlo simulations that showed how the success rate (of finding a feasible solution) and the path length change as a function of available computation time. The approach is well suited for problems involving only a few EZs (e.g., eight or less) where it can produce fairly reliable solutions (e.g., 98.8\% success rate with a computation budget of 5 seconds). For scenarios involving more EZs (e.g., between 16 and 24), especially in difficult arrangements, the algorithm will like require significantly more than 160 seconds of computation time to achieve similar success rates.

Future work may consider using the proposed approach with different $R_{\rm max}$ values to represent varying degrees of risk. This would allow selection of paths based on a trade-off of acceptable risk and path length, as well as computation of path sensitivity to risk threshold.
\bibliographystyle{IEEEtran}
\bibliography{refs} 

\begin{thebibliography}{22}
\providecommand{\natexlab}[1]{#1}
\providecommand{\url}[1]{\texttt{#1}}
\providecommand{\urlprefix}{URL }
\expandafter\ifx\csname urlstyle\endcsname\relax
  \providecommand{\doi}[1]{doi:\discretionary{}{}{}#1}\else
  \providecommand{\doi}{doi:\discretionary{}{}{}\begingroup \urlstyle{rm}\Url}\fi

\bibitem[{Boardman et~al.(2019)Boardman, Harden, and Martínez}]{boardman2018improved}
Boardman, B., Harden, T., and Martínez, S. (2019).
\newblock Improved performance of asymptotically optimal rapidly exploring random trees.
\newblock \emph{J. Dyn. Sys., Meas., Control.}, 141(1), 011002.

\bibitem[{Dillon et~al.(2023)Dillon, Zollars, Weintraub, and Von~Moll}]{dillon2023optimal}
Dillon, P.M., Zollars, M.D., Weintraub, I.E., and Von~Moll, A. (2023).
\newblock Optimal trajectories for aircraft avoidance of multiple weapon engagement zones.
\newblock \emph{J. Aerosp. Inf. Syst.}, 20(8), 520--525.

\bibitem[{Dubins(1957)}]{dubins1957curves}
Dubins, L.E. (1957).
\newblock On curves of minimal length with a constraint on average curvature, and with prescribed initial and terminal positions and tangents.
\newblock \emph{Amer. J. Math.}, 79(3), 497--516.

\bibitem[{Gammell et~al.(2020)Gammell, Barfoot, and Srinivasa}]{gammell2020batch}
Gammell, J.D., Barfoot, T.D., and Srinivasa, S.S. (2020).
\newblock Batch informed trees ({BIT*}): Informed asymptotically optimal anytime search.
\newblock \emph{The Int. J. Robot. Res.}, 39(5), 543--567.

\bibitem[{Gammell et~al.(2014)Gammell, Srinivasa, and Barfoot}]{gammell2014informed}
Gammell, J.D., Srinivasa, S.S., and Barfoot, T.D. (2014).
\newblock Informed {RRT}: Optimal sampling-based path planning focused via direct sampling of an admissible ellipsoidal heuristic.
\newblock In \emph{2014 IEEE/RSJ Int. Conf. Intell. Robots Syst.}, 2997--3004.

\bibitem[{Herrmann(1996)}]{herrmann1996air}
Herrmann, J. (1996).
\newblock Air-to-air missile engagement analysis using the {USAF} trajectory analysis program ({TRAP}).
\newblock In \emph{Flight Simul. Technol. Conf.}, 148--158.

\bibitem[{Janson et~al.(2015)Janson, Schmerling, Clark, and Pavone}]{janson2015fast}
Janson, L., Schmerling, E., Clark, A., and Pavone, M. (2015).
\newblock Fast marching tree: A fast marching sampling-based method for optimal motion planning in many dimensions.
\newblock \emph{The Int. J. Robot. Res.}, 34(7), 883--921.

\bibitem[{Karaman and Frazzoli(2011)}]{karaman2011sampling}
Karaman, S. and Frazzoli, E. (2011).
\newblock Sampling-based algorithms for optimal motion planning.
\newblock \emph{The Int. J. Robot. Res.}, 30(7), 846--894.

\bibitem[{LaValle and Kuffner~Jr(2001)}]{lavalle2001randomized}
LaValle, S.M. and Kuffner~Jr, J.J. (2001).
\newblock Randomized kinodynamic planning.
\newblock \emph{The Int. J. Robot. Res.}, 20(5), 378--400.

\bibitem[{Lozano-Perez(1983)}]{lozano1983spatial}
Lozano-Perez, T. (1983).
\newblock Spatial planning: A configuration space approach.
\newblock \emph{IEEE Trans. on Computers}, C-32(2), 108--120.

\bibitem[{O'Donohue(2019)}]{odonohu2019joint}
O'Donohue, D.J. (2019).
\newblock {Joint Air Operations (Joint Publication 3-30)}.
\newblock Technical report, Joint Chiefs of Staff, United States of America.

\bibitem[{Otte and Frazzoli(2016)}]{otte2016rrtx}
Otte, M. and Frazzoli, E. (2016).
\newblock {RRT$^X$}: Asymptotically optimal single-query sampling-based motion planning with quick replanning.
\newblock \emph{The Int. J. Robot. Res.}, 35(7), 797--822.

\bibitem[{Scott(2017)}]{scott2017countering}
Scott, K.D. (2017).
\newblock {Countering Air and Missile Threats (Joint Publication 3-01)}.
\newblock Technical report, Joint Chiefs of Staff, United States of America.

\bibitem[{Sucan et~al.(2012)Sucan, Moll, and Kavraki}]{sucan2012open}
Sucan, I.A., Moll, M., and Kavraki, L.E. (2012).
\newblock The open motion planning library.
\newblock \emph{IEEE Robot. \& Automat. Mag.}, 19(4), 72--82.

\bibitem[{Von~Moll and Weintraub(2023)}]{von2023basic}
Von~Moll, A. and Weintraub, I.E. (2023).
\newblock Basic engagement zones.
\newblock \emph{arXiv preprint arXiv:2311.06165}.

\bibitem[{Weintraub et~al.(2022)Weintraub, Von~Moll, Carrizales, Hanlon, and Fuchs}]{weintraub2022optimal}
Weintraub, I.E., Von~Moll, A., Carrizales, C.A., Hanlon, N., and Fuchs, Z.E. (2022).
\newblock An optimal engagement zone avoidance scenario in {2-D}.
\newblock In \emph{AIAA SciTech 2022 Forum}, 1587.

\bibitem[{Xu et~al.(2020)Xu, Xu, and Yin}]{xu2020optimized}
Xu, C., Xu, M., and Yin, C. (2020).
\newblock Optimized multi-{UAV} cooperative path planning under the complex confrontation environment.
\newblock \emph{Comput. Commun.}, 162, 196--203.

\bibitem[{Yan et~al.(2020)Yan, Kuang, Zhu, and Yuan}]{yan2020reachability}
Yan, X., Kuang, M., Zhu, J., and Yuan, X. (2020).
\newblock Reachability-based cooperative strategy for intercepting a highly maneuvering target using inferior missiles.
\newblock \emph{Aerosp. Sci. Technol.}, 106, 106057.

\bibitem[{Zarchan(1999)}]{zarchan1999ballistic}
Zarchan, P. (1999).
\newblock Ballistic missile defense guidance and control issues.
\newblock \emph{Sci. \& Global Secur.}, 8(1), 99--124.

\bibitem[{Zarchan(2019)}]{zarchan2012tactical}
Zarchan, P. (2019).
\newblock \emph{Tactical and Strategic Missile Guidance}.
\newblock American Institute of Aeronautics and Astronautics, Inc., 7th edition.

\bibitem[{Zhang et~al.(2023)Zhang, Jiang, Wu, and Zhu}]{zhang2023efficient}
Zhang, Z., Jiang, J., Wu, J., and Zhu, X. (2023).
\newblock Efficient and optimal penetration path planning for stealth unmanned aerial vehicle using minimal radar cross-section tactics and modified {A-Star} algorithm.
\newblock \emph{ISA Trans.}, 134, 42--57.

\bibitem[{Zhang et~al.(2020)Zhang, Wu, Dai, and He}]{zhang2020novel}
Zhang, Z., Wu, J., Dai, J., and He, C. (2020).
\newblock A novel real-time penetration path planning algorithm for stealth {UAV} in {3D} complex dynamic environment.
\newblock \emph{IEEE Access}, 8, 122757--122771.

\end{thebibliography}

\end{document}